\documentclass[11pt]{article}%
\usepackage{amsmath,amssymb}\usepackage{epsf}
\newtheorem{thm}{Theorem}
\newtheorem{co}[thm]{Corollary}
\newtheorem{lem}[thm]{Lemma}

\newtheorem{assumption}[thm]{Assumption}

\newtheorem{pr}[thm]{Proposition}

\newtheorem{definition}[thm]{Definition}

\newtheorem{example}[thm]{Example}

\newtheorem{remark}[thm]{Remark}

\newtheorem{protocol}[thm]{Protocol}

\newcommand{\Q}{\mathbb{Q}}
\newcommand{\N}{\mathbb{N}}
\newcommand{\Z}{\mathbb{Z}}

\newcommand{\Section}[1]{\section{#1}\setcounter{equation}{0}}

\newcommand{\openbox}{\leavevmode
  \hbox to.77778em{%
    \hfil\vrule
  \vbox to.675em{\hrule width.6em\vfil\hrule}%
  \vrule\hfil}}
\newcommand{\proofname}{Proof}

\newenvironment{proof}[1][\proofname]{\par\normalfont
  \trivlist\item[\hskip\labelsep\itshape #1:]\ignorespaces
  }{\hspace*{1cm}\hspace*{\fill}\openbox \medskip\endtrivlist}

\title{A New Family of Almost Identities}
\date{\today}%

\author{G\'erard Maze \\
  {\small {\em e-mail:\/} gmaze@math.unizh.ch \vspace{-2mm} }\\
  {\small Mathematics Institute\vspace{-2mm} }\\
  {\small University of Zurich\vspace{-2mm} }\\
  {\small CH-8057 Zurich,  Switzerland }\vspace{3mm}
  \and  
   Lorenz Minder \\
  {\small {\em e-mail:\/} lorenz.minder@epfl.ch \vspace{-2mm} }\\
  {\small Laboratoire de Math\'ematiques Algorithmiques\vspace{-2mm}}\\
  {\small Ecole Polytechnique F\'ed\'erale de Lausanne\vspace{-2mm}}\\
  {\small CH-1015 Lausanne, Switzerland\vspace{-2mm}}
 }
  
\begin{document} \maketitle
\thispagestyle{empty}



\Section{Introduction}                  \label{Sect:S1}

It is well-known that a class of ``almost integers'' can be
found using the theory of modular functions, and a few 
spectacular examples are given by Ramanujan \cite{raman}.  They can be
generated using some amazing properties of the $j$-function. Some of the
numbers which are close approximations of integers are $\exp(\pi
\sqrt{163})$ (sometimes known as Ramanujan's constant), $\exp(\pi
\sqrt{37})$ and $\exp(\pi \sqrt{58})$. These irrationals come close to an
integer as follows:
\begin{eqnarray*}
\exp(\pi \sqrt{37}) & = & 199148648 - 0.219... \; \cdot 10^{-4}\\
\exp(\pi \sqrt{58}) & = & 24591257752 - 0.177...\; \cdot 10^{-6}\\
\exp(\pi \sqrt{163}) & = & 262537412640768744 - 0.749...\; \cdot 10^{-12}
\end{eqnarray*}
Another surprising result comes from the average length of a segment
in an isosceles right triangle with catheti of unit length. If $l$ is this average
length, then 
$$
l=\frac{1}{30}\left( 2+4\sqrt{2}+(4+\sqrt{2})\sinh^{-1}(1) \right)
  =0.4142933026... = 
  (\sqrt{2}-1) - 0.8... \; \cdot 10^{-4}.
$$ 
Such astonishing non-equalities are usually called almost
identities or non-identities. Many examples of such unexpected
behaviour are known \cite{weis}. The four examples above are however
different in essence: the first three come from a deep property of a
complex mathematical object (the $j$-function) and the last has a good
chance to be a genuine arithmetical coincidence.

A natural question that comes to mind in presence of such a
non-identity is therefore whether or not the phenomenon is purely
coincidental, or comes from a more subtle process. For instance, in
the equation 
$$
e^{\pi} -\pi = 19.999099979...\;,
$$
it is not clear at all whether the almost identity pops up from a deep
connection between $e$ and $\pi$ or just because the expression
\textit{happens} to be close to $20$.

Recently, J.M. Borwein and P.B. Borwein discovered several families of
almost identities \cite{bor}, leading to a systematic study of such
phenomena. These were based on mathematical concepts that lead  to
clear explanations. Among the non-identities studied by these authors,
let us mention the following striking example:
$$
 \sum_{k=-\infty}^{\infty} \frac{1}{10^{(k/100)^2}} \cong  100
\sqrt{\frac{\pi}{\ln(10)}},
$$
correct to at least 18,000 digits. In this situation, the almost
identity is \textit{not} a coincidence. From the same viewpoint, let us
mention as well the sequence 
$$
h_n=\frac{n!}{2(\ln(2))^{n+1}},
$$
for $1\leqslant n \leqslant 17$, discovered by D. Hickerson.
These numbers are close to integers due to the fact that the
above quotient is the dominant term in an infinite series whose
sum is the number of possible outcomes of a race between $n$
people (where ties are allowed). See \cite{weis} for the exact
expression of
these numbers. Here, once again, no coincidence.\\

While we were studying the function
$$
f(x) = \sum_{k=1}^{\infty} \frac{1}{1+2^kx}\;, \;\;\; x \in \; ]0,1],
$$ 
that appears in the analysis of the complexity of the binary gcd
algorithm, we came to find a new family of almost identities. Let us
define the real numbers $u_n$ as follows:
$$ 
u_n := \ln(2) \cdot \sum_{k= -\infty}^{\infty}
\frac{1}{\left(2^{k/2}+2^{-k/2}\right)^n} \;\; , \;\;\; n \in
\mathbb{N}-\{0\}.
$$ 
The following equalities show the very strange behaviour of the
almost identities generated by the sequence $\{u_n\}$.
$$
\begin{array}{ccc}
u_1 & = & \displaystyle \pi \; + \; 0.53... \cdot 10^{-11}\vspace{2mm}\\ 
u_2 & = & 1 \; + \; 0.48... \cdot 10^{-10}\vspace{2mm}\\ 
u_3 & = & \displaystyle \frac{\pi}{2^3} \; + \; 0.22... \cdot
10^{-9}\vspace{2mm}\\  
u_4 & = & \displaystyle \frac{1}{6} \; +  \; 0.67... \cdot
10^{-9}\vspace{2mm} \\  
u_5 & = & \displaystyle \frac{3\pi}{2^7} \; + \; 0.15...  \cdot
10^{-8}\vspace{2mm} \\  
u_6 & = & \displaystyle \frac{1}{30} \; + \; 0.29... \cdot
10^{-8}\vspace{2mm}\\  
\vdots & \vdots &
\;\;\;\;\;\;\;\;\;\;\;\;\;\vdots \\
\end{array}
$$

This article presents an explanation of this phenomenon and
sheds light on the relation between $u_n$ and $u_{n+2}$. We first
study the cases with $n=1$ and $n=2$ by using the theory of Mellin
transforms. From there, we exhibit the recurrence relation
$$
u_n= \left(\frac{1}{4} \cdot \frac{n-2}{n-1}\right) u_{n-2} + r_n
$$
and give the explicit values of $r_n$ satisfying $0< r_n \leqslant r_{10} =
0.7227399... \cdot 10^{-8}$, $\forall n \in \N$. We also present a
generalization of the phenomenon, leading to, e.g., the almost-identity
$$
\ln(4) \cdot \sum_{k=-\infty}^{\infty}
\frac{1}{2^{-k}+2^{k}} = \pi \; + \;   0.82...
\cdot 10^{-5}.
$$

In this article, we will use the notation $f(x) \sim_a g(x)$ for
equivalent functions in a neighbourhood of $a$ and $\log_2 x$ for
the logarithm in base 2 of $x$. Also, the set $\N$ is considered to
contain the integer $0$ in the sequel.

\Section{The cases $n=1$ and $n=2$}     \label{Sect:S2}

The first two cases in our list are
$$
u_1 = \ln(2) \cdot \sum_{k=-\infty}^{\infty} \frac{1}{ 2^{-k/2}+2^{k/2} }
= \pi \; +\;  0.53... \cdot 10^{-11}
$$
and 
$$
u_2 = \ln(2) \cdot \sum_{k=-\infty}^{\infty} \frac{1}{\left( 2^{-k/2}+2^{k/2}
  \right)^2} 
= 1 \; + \; 0.48... \cdot 10^{-10}.
$$
In the next section, we will see that the expression of $u_n$, $n>2$, can be
explicited based on these first two almost identities. We therefore begin
our study by these cases. Let us define the complex functions $g_1$ and $g_2$
as
$$
g_1(x) = -2 \cdot \left(\arctan (\sqrt{x}) -\pi/2\right)
\; \mbox{ and } \; g_2(x) = \frac{1}{1+x}\; , \;\; \Re x >0
$$
as well as the functions $G_1$ and $G_2$ defined as
$$
G_n(x) = \sum_{k=1}^{\infty} g_n(2^k x) \; , \;\;\Re x>0 \;\;,\;\; n=1,2.
$$ 
The convergence of $G_1$ is justified by the fact that in a
neighbourhood of $+ \infty$ we have
$$
\arctan t - \pi /2 =-\int_{t} ^{\infty} \frac{1}{1+v^2} dv = -\int_{t}
^{\infty} \left(\frac{1}{v^2} - \frac{1}{v^4} + \frac{1}{v^6} + \cdots \right)
dv = O(1/t). 
$$
The following equalities are justified because $G_1$ and $G_2$
converge uniformly on compact subsets of their domains, and therefore
the derivative can be interchanged with the sum. Therefore,
\begin{eqnarray}
\lim_{m \rightarrow \infty} \left. \frac{d}{du}
   \left[ G_n(2^{-u}) \right] \right|_{u=m} & = &  
\lim_{m \rightarrow \infty} \sum_{k=1}^{\infty}
   \left. \frac{d}{du} \left[ g_n(2^{k-u})\right]
   \right|_{u=m} \nonumber \\
 & = &  \lim_{m \rightarrow \infty}  \ln(2) \cdot 
\sum_{k=1}^{\infty} \frac{\left(2^{(k-m)}\right)^{n/2}}{\left(
  1+2^{(k-m)}\right)^n} \nonumber \\
 & = & \lim_{m \rightarrow \infty}  \ln(2) \cdot 
\sum_{k=1}^{\infty} \frac{1}{\left( 2^{-(k-m)/2}+2^{(k-m)/2}
  \right)^n} \nonumber \\
& = & u_n,\label{eq5} 
\end{eqnarray}
where the limit is understood with $m \in \N$. The game plan is then
to express the functions $G_1$ and $G_2$ in a completely different
manner in order to compute these limits. The keystone of this process
is the Mellin transform \cite{flajolet}. Recall that the Mellin
transform of a locally Lebesgue integrable function $f(x)$ over
$]0,\infty[$ is the function
$$
f^*(s) = \int_{0}^{\infty} f(x) x^{s-1} dx.
$$
The conditions $f(x) \sim_0 O(x^u)$ and $f(x) \sim_{\infty} O(x^v)$,
with $u>v$ guarantee that $f^*(s)$ exists in the strip $-u < \Re s
<-v$. Mellin's inversion formula \cite[p.13]{flajolet} states that if 
$f$ is continuous and $c \in \; ]-u,-v[$, then 
$$
f(x) = \frac{1}{2\pi i} \int_{c-i\infty}^{c+i\infty} f^*(s) x^{-s} ds,
$$
and in a neighbourhood of $0$, we have 
$$
f(x) = \sum_{\Re s_l < c} \mbox{Res}(f^*(s) x^{-s} ,s_l),
$$
where the summation is over the poles $s_l$ of the function
$f^*(s) x^{-s}$ whose real part is strictly smaller than $c$.

Let $g(x)$ be a locally Lebesgue integrable function over
$]0,\infty[$, $f(x) = \sum_{k=1}^{\infty} g(2^kx)$, and suppose
that the convergence is uniform in $]0,\infty[$. Then
\begin{eqnarray}
f^*(s) & = & \int_{0}^{\infty} \sum_{k=1}^{\infty} g(2^kx) x^{s-1}
 dx \nonumber \\
 & = & \sum_{k=1}^{\infty}  \int_{0}^{\infty} g(y) y^{s-1} 2^{-ks} dy
 \nonumber \\
 & = & \frac{g^*(s)}{2^s-1} \label{eq1}.
\end{eqnarray}

\begin{pr} For $x>0$, we have
$$
G_1(x)=-\frac{\pi}{2}-\pi\log_2(x) +\sqrt{x}S_1(x)  -  
\sum_{k=1}^{\infty} \frac{\sin(2k\pi\log_2(x))}{k \cdot \cosh
  \left(2k\pi^2/\ln(2)\right)}
$$
where $S_1(x)$ is a power series in $x$, which converges in $[0,1[$.
\end{pr}

\begin{proof}
As announced earlier, the idea is to use Mellin transforms in a back
and forth process to reveal another expression of $G_1$. Using
(\ref{eq1}), we can write
\begin{equation}\label{eq2}
G_1^*(s) =  \frac{g_1^*(s)}{2^s-1}.
\end{equation}
In order to compute $g_1^*$, recall that in a neighbourhood of $+
\infty$ we have $ \arctan t - \pi /2 = O(1/t).  $ So, we can perform
an integration by parts, as long as $\Re s \in \;]0,1/2[$:
\begin{eqnarray*}
g_1 ^* (s) & = & - 2  \int_{0}^{\infty}  \left(\arctan
(\sqrt{x}) -\pi/2\right) x^{s-1} dx\\
& = &  -2 \cdot \left[ \left(\arctan (\sqrt{x}) - \pi /2
  \right)\cdot \left. \frac{x^s}{s} \right|_{0}^{\infty} - \frac{1}{2s}
  \int_{0}^{\infty} \frac{x^{s-1/2}}{1+x} dx \right]\\
 & = & \frac{1}{s} \int_{0}^{\infty} \frac{x^{s-1/2}}{1+x} dx\\
 & = & \frac{\pi}{s \cos \pi s}.
\end{eqnarray*}
The last equality comes from the relation
$$
\int_{0}^{\infty} \frac{x^{s-1}}{1+x} dx = \frac{\pi}{\sin \pi s}.
$$
Using Mellin's inversion formula with $c=1/4$ and (\ref{eq2}), we get 
\begin{eqnarray*}
G_1(x) & = & \frac{1}{2 \pi i} \int_{1/4-i\infty}^{1/4+i\infty}
\left(\frac{\pi}{s \cos \pi s}\right) \frac{x^{-s}}{2^s-1} ds \\ &
= & \sum_{\Re s_l < 1/4}
\mbox{Res}\left(\left(\frac{\pi}{s \cos \pi s}\right)
\frac{x^{-s}}{2^s-1} ,s_l\right).
\end{eqnarray*}
The poles of the function $\left(\frac{\pi}{s \cos \pi s}\right)
\frac{x^{-s}}{2^s-1} $ can be partitioned as follows:
\begin{itemize}
\item[{\it i)}] $s=0$ is a pole of order two, 
\item[{\it ii)}] the real simple poles $-1/2 + k$, $k \in \Z$,
\item[{\it iii)}] the imaginary simple poles $2k \pi i / \ln(2)$, $k \in \Z
  \setminus \{0\}$.
\end{itemize}
The residues are then
$$
\begin{array}{cl}
-\pi \log_2(x) - \frac{\pi}{2} & \mbox{ at } s=0, \vspace{2mm}\\
-\frac{(-2)^{k+2}}{(1+2k)(2^{k+1}-\sqrt{2})} \sqrt{x} x^k & \mbox{ at }
s=-1/2-k\; , \; k \in \N, \vspace{2mm}\\
\frac{1}{2i} \cdot \frac{\exp(-2k\pi i \log_2(x))}{k \cdot \cosh
  \left(2k\pi^2/\ln(2)\right)} & \mbox{ at }
s=2k\pi i/\ln(2)\; , \; k \in \Z \setminus \{0\} ,\vspace{2mm}\\
\end{array}
$$
and the above sum becomes
$$
G_1(x)= -\frac{\pi}{2} - \pi\log_2(x) + \sum_{k=0}^{\infty}
\frac{(-2)^{k+2}}{(1+2k)(-2^{k+1}+\sqrt{2})} \sqrt{x} x^k -  
\sum_{k=1}^{\infty} \frac{\sin(2k\pi\log_2(x))}{k \cdot \cosh
  \left(2k\pi^2/\ln(2)\right)}
$$
which proves the proposition.
\end{proof}

\begin{co}\label{co1}
$
u_1= \pi + \sum_{k=1}^{\infty} \frac{2 \pi}{\cosh
   \left(2k\pi^2/\ln(2)\right)}.
$
\end{co}

\begin{proof}
Based on (\ref{eq5}), we have
\begin{eqnarray*}
u_1 & = & \lim_{m \rightarrow \infty} \left. \frac{d}{du} \left[
   G_1(2^{-u}) \right] \right|_{u=m} \\ 
& = & \pi + \lim_{u \rightarrow \infty}
   \left[e^{-u/2}S_1(e^{-u}) \right]' + \sum_{k=1}^{\infty} \frac{2
   \pi}{\cosh \left(2k\pi^2/\ln(2)\right)}
\end{eqnarray*}
and the last limit being equal to zero, the corollary is proven.
\end{proof}

The case $n=1$ is then settled since the sum on the right-hand side of
the equality of Corollary \ref{co1} is
in fact small:
$$
u_1-\pi=\sum_{k=1}^{\infty} \frac{2 \pi}{\cosh \left(2k\pi^2/\ln(2)\right)}
= 0.538914478... \cdot 10^{-11}.
$$

\begin{pr} For $x>0$, we have
$$
G_2(x)= - \frac{1}{2} - \log_2(x) + S_2(x)
- \frac{2 \pi}{\ln(2)} \sum_{k=1}^{\infty} \frac{\sin(2k \pi
  \log_2(x))}{\sinh\left(2k\pi^2/\ln(2)\right)}
$$
where $S_2(x)$ is a power series in $x$, converging in $[0,1[$ such that
$S_2(x)=0$.
\end{pr}

\begin{proof}
The proof follows the same lines as in the first case. First,
$$
g_2^*(s) =  \int_{0}^{\infty} \frac{x^{s-1}}{1+x} dx = \frac{\pi}{\sin
  \pi s},
$$
and thus, once again based on (\ref{eq1}) and (\ref{eq2}), we have
\begin{eqnarray*}
G_2(x) & = & \int_{1/2-i \infty}^{1/2+i \infty} G_2^*(s) x^{-s} ds\\ 
&= &  \int_{1/2-i\infty}^{1/2+i\infty}\left(\frac{\pi}{\sin \pi
  s}\right)  \frac{x^{-s}}{2^s-1} ds \\ 
& = & \sum_{\Re s_l < 1/2}
 \mbox{Res}\left( \left(\frac{\pi}{\sin \pi
  s}\right)  \frac{x^{-s}}{2^s-1}    ,s_l \right).
\end{eqnarray*}

\noindent
The poles of the function can be partitioned as follows:
\begin{itemize}
\item[{\it i)}] $s=0$ is a pole of order two, 
\item[{\it ii)}] the real simple poles $ k$, $k \in \Z\setminus \{0\}$,
\item[{\it iii)}] the imaginary simple poles $2k \pi i / \ln(2)$, $k \in \Z
  \setminus \{0\}$.
\end{itemize}
\noindent
The residues are then
$$
\begin{array}{cl}
- \log_2(x) - \frac{1}{2} & \mbox{ at } s=0, \vspace{2mm}\\
-\frac{(-2)^k}{2^k-1}x^k & \mbox{ at } s=-k\; , \; k =1,2,3,...\; ,
\vspace{2mm}\\ 
\frac{\pi}{i} \cdot \frac{\exp(-2k\pi i \log_2(x))}{\ln(2) \cdot \sinh
  \left(2k\pi^2/\ln(2)\right)}  & \mbox{ at }
s=2k\pi i/\ln(2)\; , \; k \in \Z \setminus \{0\}.\vspace{2mm}\\
\end{array}
$$
The new expression of $G_2$ is therefore
$$
G_2(x) = - \frac{1}{2} - \log_2(x) - \sum_{k=1}^{\infty}
\frac{(-2)^k}{2^k-1}x^k
- \frac{2 \pi}{\ln(2)} \sum_{k=1}^{\infty} \frac{\sin(2k \pi
  \log_2(x))}{\sinh\left(2k\pi^2/\ln(2)\right)} .
$$
This concludes the proof.
\end{proof}

\begin{co}\label{co2}
$
u_2= 1 + \frac{2 \pi}{\ln(2)} \sum_{k=1}^{\infty} \frac{2k
   \pi}{\sinh \left(2k\pi^2/\ln(2)\right)}.
$
\end{co}

\begin{proof}
We use here the same trick as in Corollary \ref{co1}:
\begin{eqnarray*}
u_2 & = & \lim_{m \rightarrow \infty} \left. \frac{d}{du} \left[
   G_2(2^{-u}) \right] \right|_{u=m} \\ 
& = & 1 + \lim_{u \rightarrow \infty}
   \left[S_2(e^{-u}) \right]' + \frac{2 \pi}{\ln(2)} \sum_{k=1}^{\infty} \frac{2k
   \pi}{\sinh \left(2k\pi^2/\ln(2)\right)}
\end{eqnarray*}
and the limit being equal to zero, the corollary is proven.
\end{proof}

Once again, this shows why the number $u_2$ is almost an integer. Indeed the sum on
the right-hand side is fairly small:
$$
u_2-1=\frac{2 \pi}{\ln(2)} \sum_{k=1}^{\infty} \frac{2k
   \pi}{\sinh \left(2k\pi^2/\ln(2)\right)} = 0.4885108992... \cdot 10^{-10}.
$$

\Section{The recurrence relation}        \label{Sect:S3}

Having found the roots of the mystery related to the non-equalities
$u_1 \neq \pi$ and $u_2 \neq 1$, we would now like to extend the
method used in the previous section to understand why $u_3,u_4,...$
are so close to ``good arithmetic numbers''.
Looking back to the cases $n=1,2$, we see that the functions $g_1$ and
$g_2$ played a crucial role. The key was the fact that they satisfy
the equalities
$$
\frac{d}{du} \left[ g_n(2^{k-u}) \right] = \frac{\ln 2}{\left(
  2^{-(k-u)/2}+2^{(k-u)/2} \right)^n}\;, \;\;n=1,2.
$$
The next lemma shows how we can extend them:
\begin{lem}
Let $n \in \N$, $n>2$, and let
\begin{eqnarray*}
I_{n,k} & = &  \int \frac{1}{\left( 2^{-(k-u)/2}+2^{(k-u)/2}
  \right)^n} \, du,\\
R_{n,k}  & = & \frac{1}{2 \ln2 \cdot (n-1)} \left(\frac{2^{(k-u)/2}}{1+2^{k-u}}
\right)^{n-2} \left(\frac{1-2^{k-u}}{1+2^{k-u}}\right).
\end{eqnarray*}
Then
$$
I_{n,k} = \left( \frac{1}{4} \cdot \frac{n-2}{n-1}\right)  I_{n-2,k} + R_{n,k}
$$
\end{lem}

The proof is left to the reader, who can simply differentiate and
check! The equality of the previous lemma can be used as
follows. For $n>2$, we have
\begin{eqnarray} 
u_n & = & \ln(2) \cdot \lim_{m \rightarrow \infty} \sum_{k=1}^{\infty}
  \frac{1}{\left( 
  2^{-(k-m)/2}+2^{(k-m)/2} \right)^n} \nonumber \\
 & = & \ln(2) \cdot   \lim_{m \rightarrow \infty} \sum_{k=1}^{\infty} 
  \frac{d}{du}\left[ I_{n,k}\right]_{u=m} \nonumber \\
 & = & \ln(2) \cdot  \lim_{m \rightarrow \infty}  \left(
  \sum_{k=1}^{\infty} \frac{d}{du} \left[\left( \frac{1}{4} \cdot
  \frac{n-2}{n-1}\right) I_{n-2,k} + R_{n,k} \right]_{u=m}\right)
  \nonumber \\  
 & = & \left( \frac{1}{4} \cdot \frac{n-2}{n-1}\right) u_{n-2} + \ln(2)
  \cdot \lim_{m \rightarrow \infty}  \left(
  \sum_{k=1}^{\infty} \frac{d}{du} \left[ R_{n,k}
  \right]_{u=m}\right).\label{eq3} 
\end{eqnarray}
Let us define, for $n>2$, 
$$
f_n(x)=  \left(\frac{\sqrt{x}}{1+x}\right)^{n-2}
  \frac{1-x}{1+x} \;,\;\; \mbox{ so that } \;\;\; \frac{1}{2 \ln2
    \cdot (n-1)} f_n(2^{k-u}) = R_{n,k}. 
$$
If
$$
F_n(x) = \sum_{k=1}^{\infty} f_n(2^kx),
$$
since this function converges uniformly and absolutely on compact subsets
of $\Re x>0$, we can interchange derivation and summation to obtain
\begin{eqnarray}\label{eq4}
\frac{1}{2 \cdot (n-1)}\cdot \lim_{m \rightarrow \infty}
\frac{d}{du} \left[ F_n(2^{-u})
  \right]_{u=m} & = & \ln(2) \cdot \lim_{m
  \rightarrow \infty}  \left( 
  \sum_{k=1}^{\infty} \frac{d}{du} \left[ R_{n,k} \right]_{u=m}\right).
\end{eqnarray}

Once again, we use Mellin transforms to find another expression for
each of the functions $F_n(x)$ in order to compute these limits.
\begin{pr}
The function $F_n$, $n>2$, can be represented as
$$
F_n(x) = 
\left\{
\begin{array}{cl}\displaystyle 
 S_n(x) - \frac{4 \pi}{\ln 2}
  \sum_{k=1}^{\infty}  c_k 
  \frac{\sin(2k\pi \log_2(x))}{\sinh(2k\pi^2/\ln2)}  & \mbox{ when $n$ is
  even,} \\ \displaystyle  
  \sqrt{x} S_n(x) - \frac{4 \pi}{\ln 2}
  \sum_{k=1}^{\infty}  b_k 
  \frac{\sin(2k\pi \log_2(x))}{\cosh(2k\pi^2/\ln2)}  & \mbox{ when $n$ is odd,}
\end{array}
\right.
$$
where $S_n(x)$ is a power series converging in $[0,1[$ such that
$S(0)=0$. The coefficients $c_k$ and $b_k$ are given by
$$
\begin{array}{ll}\displaystyle 
c_k=  \frac{\prod_{j=0}^{l-2} (j^2+4\pi^2k^2/\ln(2)^2)}{(2l-2)!} & \mbox{
  when } n= 2l,\; l>1, \\ \displaystyle   
b_k = \frac{2 \pi k \prod_{j=0}^{l-2}
  ((j+1/2)^2+4\pi^2k^2/\ln(2)^2)}{\ln(2)(2l-1)!}  & \mbox{ when } n=
  2l+1, \; l>1. 
\end{array}
$$
\end{pr}

\begin{proof}
First,
\begin{eqnarray*}
f_n^*(s) & = & \int_{0}^{\infty} f_n(x) x^{s-1} dx\\
 & = & \int_{0}^{\infty} \frac{x^{n/2+s-2}}{(1+x)^{n-2}}
 \frac{1-x}{1+x} dx\\
 & = &  \int_{0}^{\infty} \frac{x^{n/2+s-2}}{(1+x)^{n-1}}dx 
 - \int_{0}^{\infty} \frac{x^{n/2+s-1}}{(1+x)^{n-1}} dx.\\
\end{eqnarray*}
This expression can be evaluated with the help of the Gamma function
$\Gamma$. Indeed, this function satisfies, see, e.g., \cite[p.47]{gamma},
$$
\frac{\Gamma(p)\Gamma(q)}{\Gamma(p+q)} = \int_{0}^{\infty}
\frac{x^{p-1}}{(1+x)^{p+q}} dx
$$
and therefore
\begin{eqnarray*}
f_n^*(s) & = & \frac{\Gamma(n/2+s-1)\Gamma(n/2-s)}{\Gamma(n-1)} -
\frac{\Gamma(n/2+s)\Gamma(n/2-s-1)}{\Gamma(n-1)} \\
& = & -2s \frac{\Gamma(n/2-1+s)\Gamma(n/2-1-s)}{\Gamma(n-1)}.
\end{eqnarray*}
We used the equality $\Gamma(n)=(n-1)\Gamma(n-1)$ in the last
step. Based on Euler's reflection formula $\Gamma(s)\Gamma(1-s) =
\pi/\sin(\pi s)$, see, e.g., \cite[p.9]{gamma}, the previous equality
leads to the following expressions, both correct for $\Re s \in
]0,1/2[$ :
$$
f_n^*(s)=
\left\{
\begin{array}{cl}\displaystyle 
 \frac{2}{(2l-2)!} \cdot \frac{\pi }{\sin \pi s} \prod_{j=0}^{l-2}
  (j^2-s^2) & \mbox{ when } n= 2l,\; l>1, \\ \displaystyle 
 \frac{2}{(2l-1)!} \cdot  \frac{-\pi s}{\cos \pi s} \prod_{j=0}^{l-2}
  \left((j+1/2)^2-s^2\right) & \mbox{ when } n= 2l+1, \; l \geqslant 1.
\end{array}
\right.
$$

The equality (\ref{eq1}) and (\ref{eq2}) lead once again to
\begin{eqnarray*}
F_n(x) & = & \int_{1/4-i \infty}^{1/4+i \infty} F_n^*(s) x^{-s} ds\\ 
& = & \sum_{\Re s_l < 1/4}
 \mbox{Res}\left( f_n^*(s)  \frac{x^{-s}}{2^s-1}    ,s_l \right).
\end{eqnarray*}

The poles of the function can be partitioned as follows:
\begin{itemize}
\item[{\it i)}] the imaginary simple poles $2k \pi i / \ln(2)$, $k \in \Z
  \setminus \{0\}$,
\item[{\it ii)}] the real simple poles:
 \begin{itemize}
   \item[] $s=j \in \Z $ with $|j|>l-2$ when $n=2l$,
   \item[] $s=j+1/2 \in \Z+1/2 $ with $|j|>l-2$ when $n=2l+1$.
\end{itemize}
\end{itemize}

Note that contrary to the cases we have seen so far, there
are no poles with multiplicity. The real simple poles will clearly
contribute to residues of the form $a_kx^k$ when $n$ is even and $a_k
\sqrt{x}x^k$ when $n$ is odd. We do not exhibit the coefficients $a_k$
since we will not need them. The imaginary simple poles lead to 
residues at $s= 2k\pi i/ \ln 2$ , $k\neq0$, which are of the form
$$
 \mbox{Res}\left( f_n^*(s)  \frac{x^{-s}}{2^s-1}, 2k\pi i/ \ln 2\right)=  
f_n^*( 2k\pi i/ \ln 2)  \frac{x^{- 2k\pi i/ \ln 2}}{\ln(2)}.
$$
The new expression of $F_n$ is therefore
$$
F_n(x) = 
\left\{
\begin{array}{cl}\displaystyle 
 \sum_{k=l-1}^{\infty} a_kx^k - \frac{4\pi}{\ln 2}
  \sum_{k=1}^{\infty}  c_k 
  \frac{\sin(2k\pi \log_2(x))}{\sinh(2k\pi^2/\ln2)}  & \mbox{ when
  } n= 2l,\; l>1, \\ \displaystyle  
  \sqrt{x} \sum_{k=l-1}^{\infty} a_kx^k - \frac{4 \pi}{\ln 2}
  \sum_{k=1}^{\infty}  b_k 
  \frac{\sin(2k\pi \log_2(x))}{\cosh(2k\pi^2/\ln2)}  & \mbox{ when } n=
  2l+1, \; l\geqslant 1, 
\end{array}
\right.
$$
where the coefficients $c_k$ and $b_k$ are given in the proposition.
This finishes the proof.
\end{proof}

\begin{co}\label{co3}
The sequence $\{u_n\}_{n \in \N}$ satisfies the following recurrence relation
$$
u_n= \left(\frac{1}{4} \cdot \frac{n-2}{n-1}\right) u_{n-2} + r_n
$$
where 
$$
r_n= 
\left\{
\begin{array}{cl}\displaystyle 
  \frac{2 \pi}{\ln(2) (n-1)} \cdot \sum_{k=1}^{\infty}  c_k 
  \frac{2k\pi}{\sinh(2k\pi^2/\ln2)}  & \mbox{ when
  } n= 2l,\; l>1, \\ \displaystyle  
  \frac{2 \pi}{\ln(2) (n-1)} \cdot  \sum_{k=1}^{\infty}  b_k 
  \frac{2k\pi}{\cosh(2k\pi^2/\ln2)}  & \mbox{ when } n= 
  2l+1, \; l\geqslant 1.
\end{array}
\right.
$$
\end{co}

\begin{proof}
Based on (\ref{eq3}), (\ref{eq4}) and the previous
proposition, we have   
$$
u_n- \left(\frac{1}{4} \cdot \frac{n-2}{n-1}\right) u_{n-2}  = 
   \frac{1}{2 \cdot (n-1)}\cdot \lim_{m \rightarrow \infty} 
   \frac{d}{du} \left[ F_n(2^{-u})\right]_{u=m}.
$$
The limit in the above expression annihilates the limit
of the power series of $F_n$ and the only contributing term in the
limit is the sinus series of $F_n$. This gives the expected expression of $r_n$.
\end{proof}

The growth of the coefficients $r_n$ is the combined effect of the
increase of the values of $c_k$ and $b_k$ and the decrease of
$(n-1)^{-1}$. As a consequence, the sequence $r_n$ is increasing for
$n \leqslant 10$ and decreasing for $n \geqslant 10$, which gives
$$
0<r_n \leqslant r_{10} = 0.7227399... \cdot 10^{-8}.
$$

We end this article by the following remark. The entire theory used
here to explain why the numbers $u_n$ are so close to elements in $\Q
\cup \pi \Q$ has nothing to do with the presence of $2$ in the
denominator of 
$$
\frac{1}{\left(2^{-k/2}+2^{k/2}\right)^n}.
$$
One could argue that any sum of the type 
$$
\ln(m) \cdot \sum_{k=-\infty}^{\infty}
\frac{1}{\left(m^{-k/2}+m^{k/2}\right)^n}
$$
has the potential to lie close to $\Q$ or $\pi\Q$ depending on the
parity of $n$. As a matter of fact, we have, for example,
\begin{eqnarray*}
\ln(4) \cdot \sum_{k=-\infty}^{\infty}
\frac{1}{2^{-k}+2^{k}} & = & \pi \; + \;   0.82...
\cdot 10^{-5},\\
\ln(9) \cdot \sum_{k=-\infty}^{\infty}
\frac{1}{3^{-k}+3^{k}} & = & \pi \; + \;   0.15...
\cdot 10^{-2},\\
\ln(4) \cdot \sum_{k=-\infty}^{\infty}
\frac{1}{\left(2^{-k}+2^{k}\right)^2} & = & 1 \; + \; 0.37...
\cdot 10^{-4}.
\end{eqnarray*}
Based on what has been shown in this article, we can say that the
``error'' term is due to the size of $\ln(m)$ (in the hyperbolic
functions of $r_n$) and the smaller it is, the smaller the error will be. In
other words, the choice $m=2$ is the best one can do in order to
maximize the resemblance with elements in $\Q \cup \pi \Q$.

\end{document}